\documentclass[12pt]{article}

\usepackage{amsfonts,epsfig}
\usepackage{epsfig}
\usepackage{amsmath}
\usepackage{amssymb}
\usepackage{hyperref}

\newtheorem{theorem}{Theorem}
\newtheorem{corollary}[theorem]{Corollary}

\def\perm{\textnormal{perm}}
\def\trace{\textnormal{trace}}
\def\sgn{\textnormal{sgn}}
\def\term{\textnormal{term}}
\def\product{\textnormal{product}}
\def\prob{\textnormal{prob}}

\def\qedspace{\null}
\def\qedblack{\qedspace                    
\hbox{\vrule height7pt width 5pt}}

\def\QED{\qedblack}

\def\BeginProofOf#1{\noindent{\textbf{Proof of #1\hskip1pt}:} }
\def\EndProof{\QED\betweenskip}
\def\betweenskip{\vskip10pt}

\title{The expected characteristic and permanental polynomials
of the random Gram matrix}

\author{Jacob G. Martin\\
\small Department of Neuroscience\\[-0.8ex]
\small Georgetown University Medical Center\\[-0.8ex]
\small Washington, DC 20007 \\
\small\tt jm733@georgetown.edu\\
\and
E. Rodney Canfield\thanks
{Research supported by NSA Mathematical Sciences Program.}\\
\small Department of Computer Science\\[-0.8ex]
\small University of Georgia\\[-0.8ex]
\small Athens, GA 30602\\
\small\tt erc@cs.uga.edu\\}

\begin{document}

\maketitle

\begin{abstract}
A $t \times n$ random matrix $A$ can be formed by sampling $n$
independent random column vectors, each containing $t$  components.
The \textit{random Gram matrix} of size $n$, $G_{n}=A^{T}A$,
contains the dot products between all pairs of column vectors
in the randomly
generated matrix $A$, and has characteristic roots coinciding
with the singular values of $A$.  Furthermore, the sequences
$\det{(G_{i})}$
and $\perm(G_{i})$ (for $i = 0, 1, \dots, n$) are
factors that comprise  the expected coefficients of the characteristic
and permanental polynomials of $G_{n}$.  We prove theorems that
relate the generating functions and recursions for the traces
of matrix powers, expected characteristic coefficients, expected
determinants $E(\det{(G_{n})})$, and
expected permanents $E(\perm(G_{n}))$
in terms of each other.  Using the derived recursions, we exhibit
the efficient computation of the expected determinant and expected
permanent of a random Gram matrix $G_{n}$, formed according to
any underlying distribution.  These theoretical results may be
used both to speed up numerical algorithms and to investigate
the numerical properties of the expected characteristic and permanental
coefficients of any matrix comprised of independently sampled
columns.
\end{abstract}

\section{Introduction}
Let $w$ be a $t$-tall vector whose components $w_{i}$ are random
variables (not necessarily independent)
\begin{equation*}
 w = \left[ 
         \begin{array}{c}
            w_{1} \\
 	    w_{2} \\
	    \vdots \\
	    w_{t} 
         \end{array} 
         \right] \; .
\end{equation*}
Next, sample $n$ independent vectors $w^{(1)}, \dots, w^{(n)}$;
creating a $t\times n$ matrix $A$.  Then, the \textit{random
Gram matrix} of size $n$,
\begin{equation*}
 G_{n} = A^{T}A \; ,
\end{equation*}
has a distribution that depends on the underlying distribution
of the random vector $w$.  (The symbol $T$ as a superscript is used to denote
transpose.)  

Some general features and convergence properties of the eigenvalues
of certain random Gram matrices were derived by Fannes and Spincemaille
\cite{Fannes}.  Fyodorov formulated correlation functions for
permanental polynomials of certain random matrices and noted
some similarities and differences between their characteristic
and permanental polynomials \cite{Fyodorov}.  Our paper presents
combinatorial theory and an efficient algorithm for calculating
$E(\det(G_{n}))$ and $E(\perm(G_{n}))$, which are
factors comprising the coefficients of the expected characteristic
and expected permanental polynomials of $G_{n}$.  

The computation of the determinant is equivalent to matrix multiplication
and is therefore contained in the complexity class P (see Chapter
16 of \cite{Burgisser97}).  Currently, the fastest asymptotic
algorithm for matrix multiplication is $O(n^{2.376})$ \cite{Coppersmith},
with a recent unpublished work \cite{Williams} claiming an improvement to $O(n^{2.3727})$.
Some researchers have suggested that group theoretic observations
imply that $O(n^{2})$ algorithms also exist \cite{Umans}.  

At the other complexity extreme, even though
the sign is the only difference between the
formula for the determinant,
\begin{equation}
\det(A) = \sum_{\sigma \in S_{n}} (-1)^{\sgn(\sigma)} \prod_{i=1}^{n}
A_{i,\sigma(i)}, 
\end{equation}
and the formula for the permanent,
\begin{equation}
\perm(A) = \sum_{\sigma \in S_{n}} \prod_{i=1}^{n}
A_{i,\sigma(i)}, 
\end{equation}
the computation of the permanent is \#P-Complete \cite{Valiant,BenDor}.
  The standard reference for properties of permanents is
Minc \cite{Minc}.  The most efficient algorithm currently
known for calculating the exact permanent has complexity $O(2^{n}n^{2})$,
due to Ryser \cite{Ryser}.  Jerrum and Sinclair provided
a fully-polynomial randomized approximation scheme for approximating
permanents of nonnegative matrices \cite{Jerrum1,Jerrum2}.   Matrix permanents have found
applications in physics for calculating Bose-Einstein corrections
\cite{Wosiek} and in quantum computing for encoding quantum circuit
amplitudes \cite{Rudolph, Broome}.  Permanental polynomials have been used
as invariants for chemical structures \cite{Cash,Cash2,Merris}.

The rest of the paper
is organized as follows: Section 2 is a statement of results,
Section 3 contains all proofs, Section 4 reports on some numerical
experiments, Section 5 points out a connection to prior work involving
the cycle index polynomial of the symmetric group, and Section 6 presents
summary and conclusions.

\section{Statement of Results}
Before stating our
results, we explain all notation.
Let $w$ be a $t$-tall vector whose components
$w_{i}$ are random variables (not necessarily independent). 
Let $A$ be a $t \times n$ matrix whose columns are
a random sample of
$n$ vectors $w^{(1)},\dots,w^{(n)}$.
Let $G_{n}=A^{T}A$;
we call $G_n$ {\it the random
Gram matrix of size} $n$, it being understood that the exact distribution
of $G_n$ depends on the underlying distribution on $t$-dimensional vectors
$w=(w_1,\dots,w_t)^T$.   Although $G_n$ is an $n\times n$ matrix, its rank is at most $t$,
and generally speaking we take the viewpoint henceforth that $n$ is
much larger than $t$.   One may even regard $t$ as fixed, and $n\rightarrow\infty$,
as we study the effect of taking larger and larger samples.  We are especially
interested in two expected values, the determinant and the permanent of $G_n$;
these are denoted $a_n,p_n$ respectively:
\begin{eqnarray*}
a_{n} &=& E(\det(G_{n}))  \\
p_{n} &=& E(\perm(G_{n})).
\end{eqnarray*}
We define $M$ to be the $t\times t$ matrix of underlying second moments,
$$
M_{ij} = E(w_{i}w_{j}),  ~~ 1\le i,j\le t,
$$
and define the infinite sequence $t_n$ as the traces of the powers of $M$:
$$
t_{n} = \trace(M^{n}).
$$
Finally, we define $c_i$, $0\le i\le t$, to be the
sign-adjusted coefficients of the characteristic
polynomial of $M$, with the familiar indexing: 
$$
\det(\lambda I - M) = c_0\lambda^t - c_1\lambda^{t-1} + \cdots + (-1)^t c_t.
$$

\begin{theorem}{\label{theorem1}}
Let $a_n$, $p_n$, $t_n$ denote $E(\det(G_{n}))$,
$E(\perm(G_{n}))$, $\trace(M^n)$, respectively,
as given above.  Then,
\begin{align}\label{eqgfdet}
\sum_{n=0}^{\infty}a_{n}\frac{x^{n}}{n!} &= \exp{\left\lbrace\frac{t_{1}x}{1}-\frac{t_{2}x^{2}}{2}+\frac{t_{3}x^{3}}{3}-
\dots\right\rbrace} .
\end{align}
and
\begin{align}\label{eqgfperm}
\sum_{n=0}^{\infty}p_{n}\frac{x^{n}}{n!} &= \exp{\left\lbrace\frac{t_{1}x}{1}+\frac{t_{2}x^{2}}{2}+\frac{t_{3}x^{3}}{3}+
\dots\right\rbrace} .
\end{align}
\end{theorem}

\noindent The generating function identities in the previous theorem lead immediately
to recursions for the sequences $a_n,p_n$ as given in the corollary:
\begin{corollary}{\label{corollary2}}
Let $a_n$, $p_n$, $t_n$ denote $E(\det(G_{n}))$,
$E(\perm(G_{n}))$, $\trace(M^n)$, respectively,
as given above.  Then, we have the recursions
\begin{align}
	a_{0} &= 1  \nonumber \\
	a_{n+1} &= \sum_{j}\binom{n}{j}(-1)^{j}j!a_{n-j}t_{j+1}
                    \label{eqdetrecursion}
\end{align}
and 
\begin{align}
	p_{0} &= 1  \nonumber \\
	p_{n+1} &= \sum_{j}\binom{n}{j}j!p_{n-j}t_{j+1} 
         \label{eqpermrecursion}
\end{align}
\end{corollary}

\noindent
The next theorem relates the expected values $E(\det(G_n)),E(\perm(G_n))$
to the coefficients $c_i$ of the characteristic polynomial for $M$.
\begin{theorem}{\label{theorem3}}  Let $G_n$, $M$, $c_n$ be respectively
the random Gram matrix of size $n$, the underlying $t\times t$
matrix of second moments, and the sign-adjusted coefficients of the
characteristic polynomial,
$$
\det(\lambda I - M) = c_0\lambda^t - c_1\lambda^{t-1} + \cdots + (-1)^tc_t.
$$
Then
$$
E(\det(G_{n})) = n! c_n.
$$
and
$$
E(\perm(G_{n})) = n! \times [x^n] \, (1-c_1x + c_2x^2 - \cdots)^{-1}.
$$
\end{theorem}

\vskip 10pt

\noindent  The last theorem concerns the expected values of the coefficients
of the characteristic and permanental polynomials of $G_n$.
\begin{theorem}{\label{theorem4}}
Let $a_n$, $p_n$ denote $E(\det(G_{n}))$, $E(\perm(G_{n}))$, respectively,
as given above.  Let $b_i,d_i$ be the sign-adjusted coefficients of,
respectively, the characteristic and permanental polynomials $G_n$:
\begin{eqnarray*}
\det (\lambda I - G_n) &=& b_0\lambda^n - b_1\lambda^{n-1} + \cdots + (-1)^nb_n \\
\perm(\lambda I - G_n) &=& d_0\lambda^n - d_1\lambda^{n-1} + \cdots + (-1)^nd_n.
\end{eqnarray*}
Then,
\begin{equation}
\label{eqcharcoefexpect}
E(b_i) = {n\choose i} a_i
\end{equation}
and
\begin{equation}
\label{eqpermcoefexpect}
E(d_i) = {n\choose i} p_i.
\end{equation}
\end{theorem}

\vskip 10pt

\noindent{\bf Remark.}  The characteristic polynomials $\det(\lambda I - AA^T)$
and $\det(\lambda I - A^TA)$ have exactly (including multiplicity) the same nonzero
roots.  With $A$ a $t\times n$ matrix, and assuming $n\ge t$, then, the latter 
characteristic polynomial has a factor of $\lambda^{n-t}$, and so $b_i=0$ for
$i>t$.  This is consistent with the fact that the $a_i$ are nonzero for at most
$0\le i\le t$.

\section{Proofs}

\BeginProofOf{Theorem \ref{theorem1}} 
The Leibniz formula for the determinant is 
\begin{equation}\label{determinant}
 \det(G_{n}) = \sum_{\sigma \in S_{n}} (-1)^{\sgn(\sigma)}
\term_{\sigma} \; ,
 \end{equation}
with 
\begin{equation*}
 \term_{\sigma} = \prod_{i=1}^{n} (G_{n})_{i,\sigma(i)}
\; , 
\end{equation*}
where $S_{n}$ is the \textit{symmetric group} on $\left\lbrace
1,2,\dots, n \right\rbrace$, and $\sgn(\sigma)$ signifies the
\textit{sign} of the permutation $\sigma$.
Similarly, for the permanent,
\begin{equation}\label{permanent}
 \perm(G_{n}) = \sum_{\sigma \in S_{n}} \term_{\sigma}
\; .
\end{equation}

Since expectation is a linear operator, $a_{n} = E(\det(G_{n}))$
may be obtained by the following strategy
\begin{equation*}
 (*) \begin{cases}
\textnormal{1.  Determine } E(\term_{\sigma}) \\
\textnormal{2.  Multiply by } (-1)^{\sgn(\sigma)} \\
\textnormal{3.  Sum over } \sigma \in S_{n} \; .                
                \end{cases}
\end{equation*}  Furthermore, $p_{n} = E(\perm(G_{n}))$
can be obtained in the same manner but omitting step $2$.  

Suppose the permutation $\sigma$ contains $k_i$ cycles of size $i$,
where $k_i\ge 0$ and $n=k_1+2k_2+\cdots$.
The cycle structure of $\sigma$ alone is enough to determine its
sign by the relation: $\sgn(\sigma)=k_2+k_4+\cdots$.  What can be said about
the expected value
$E(\term_{\sigma})$, given only the cycle structure of
$\sigma$ ?  We claim that, like the sign, the latter expected value
is determined completely by the cycle structure, as given in
the relation
\begin{equation}
\label{eqexpectterm}
E(\term_{\sigma}) = \prod_{i\ge 1} (t_i)^{k_i}.
\end{equation}
Indeed, if $(i_{1}, i_{2}, \dots, i_{l})$ is a cycle
of $\sigma$ having length $l$, then the quantity
$\product_C$, defined by
$$
 \product_C = (G_{n})_{i_{1},i_{2}}(G_{n})_{i_{2},i_{3}}\dots (G_{n})_{i_{l},i_{1}}
\; \; ,
$$
is a \textit{subproduct} of $\term_{\sigma}$.  Moreover, it is seen
that the various subproducts
associated with the different cycles comprising $\sigma$ have no rows or columns
of the matrix $G_n$ in common.  These subproducts are consequently independent,
and we have 
\begin{equation*}
E(\term_{\sigma}) = \prod_{C:C\textnormal{ is a cycle of
} \sigma} E(\product_{C}) \;.
\end{equation*}
The entry $(G_n)_{ij}$ of the Gram matrix is the dot product
$w^{(i)} \cdot w^{(j)}$  of columns in the sample matrix $A$, and so
\begin{equation*}
 \product_{C} = (w^{(i_{1})} \cdot w^{(i_{2})}) \times
(w^{(i_{2})} \cdot w^{(i_{3})}) \times \cdots \times (w^{(i_{l})} \cdot
w^{(i_{1})}) \; \; .
\end{equation*}
From this we observe that the expectation $E(\product_{C})$
depends only on the length of the cycle
$C$, and not on the particular
columns of $A$ which are involved.  That the  
common value of $E(\textnormal{product}_{C})$ over all cycles
$C$ of length $l$ is equal to $t_{l}$, the trace of the power
$M^{\ell}$  is seen as follows
\begin{align*}
      ~ &~~~~ E\left[ \left(w^{(1)} \cdot w^{(2)}\right) \left(w^{(2)}
\cdot w^{(3)}\right) \dots \left(w^{(n)} \cdot w^{(1)}\right)\right] \\
  &= E \biggl[
     \left(w_{1}^{(1)}w_{1}^{(2)} + \dots +w_{t}^{(1)}w_{t}^{(2)}\right)
~\times~
     \left(w_{1}^{(2)}w_{1}^{(3)} + \dots +w_{t}^{(2)}w_{t}^{(3)}\right)
~\times~ \cdots \\
    &~~~~~ \cdots 
~\times~  \left(w_{1}^{(n)}w_{1}^{(1)} + \dots+w_{t}^{(n)}w_{t}^{(1)}\right)\biggr]\\
  &= E \left[ \sum_{\left(i_{1}, \dots , i_{n}\right) \in \left\lbrace
1,\dots ,t \right\rbrace^{n}} w_{i_{1}}^{(1)}w_{i_{1}}^{(2)}
w_{i_{2}}^{(2)}w_{i_{2}}^{(3)} \dots w_{i_{n}}^{(n)}w_{i_{n}}^{(1)}\right]\\
	&= \sum_{\left(i_{1}, \dots , i_{n}\right) \in \left\lbrace
1,\dots ,t \right\rbrace^{n}} E\left(w_{i_{1}}^{(1)}w_{i_{n}}^{(1)}\right)
E\left(w_{i_{1}}^{(2)}w_{i_{2}}^{(2)}\right) 
E\left(w_{i_{2}}^{(3)}w_{i_{3}}^{(3)}\right) \dots
E\left(w_{i_{n-1}}^{(n)}w_{i_{n}}^{(n)}\right) \\
 	&= \sum_{\left(i_{1}, \dots , i_{n}\right) \in \left\lbrace
1,\dots ,t \right\rbrace^{n}} M_{i_{1}i_{n}}M_{i_{1}i_{2}}M_{i_{2}i_{3}}\dots
M_{i_{n-1}i_{n}} \\
 	&= \trace(M^{n}) \\
 	&= t_n \; \; .
\end{align*}
Thus the claim (\ref{eqexpectterm}) is justified.  

Continuing the proof, we introduce the
sign of the permutation to obtain
\begin{equation}
\label{eqsignedexpect}
(-1)^{\sgn(\sigma)} E(\term_{\sigma}) = (-1)^{k_2+k_4+\cdots} ~~~ \prod_{i\ge 1} (t_i)^{k_i}.
\end{equation}
We are now in position to
carry out the three-step strategy (*) proposed above.  The number of permutations $\sigma$ which have
a given cycle structure $(k_1,k_2,\dots)$ is
\begin{equation*}
\frac{n!}{1^{k_{1}}2^{k_{2}}\dots k_{1}!k_{2}! \dots} \; \;, ~~~~ n=k_1+2k_2+\cdots \; .
\end{equation*}
If we multiply the right side of (\ref{eqsignedexpect}) by the latter multiplicity and
by $x^n/n!$ -- note the resulting cancellation of $n!$ -- and then sum over all sequences
$(k_1,k_2,\dots)$ of nonnegative integers which are zero from some point on, we obtain the
desired exponential generating function.  Hence,
\begin{align*}
 &~~~~ \sum_{n=0}^{\infty} a_n \frac{x^n}{n!} \\
 &= \sum_{(k_1,k_2,\dots)} (-1)^{k_2+k_4+\cdots} x^{k_1+2k_2+\cdots}
\prod_{i\ge 1} \frac{t_i^{k_i}}{i^{k_i}k_i!} \\
 &= \prod_{i\ge 1} \sum_{k_i=0}^{\infty} \frac{((-1)^{i-1}t_ix^i/i)^{k_i}}{k_i!} \\
 &= \exp(\frac{t_1x}{1} - \frac{t_2x^2}{2} + \cdots),
\end{align*}
as was to be shown.  This
completes the proof of the first part of Theorem~\ref{theorem1}.
The second part of the theorem, equation (\ref{eqgfperm}) giving
the exponential generating function of the
sequence of permanents $p_n$, is proven in a similar manner.
\EndProof

\BeginProofOf{Corollary \ref{corollary2}}
These are proven in the standard manner by comparing the coefficients of $x^n$
on both sides of the identities obtained from (\ref{eqgfdet})
and (\ref{eqgfperm}) by
differentiating with respect to $x$.
\EndProof

\noindent In the next proof of Theorem \ref{theorem3} we use the identity
\begin{equation}
\label{eqexptrace}
\det(\exp(B)) = \exp(\trace(B)),
\end{equation}
valid for any complex square matrix $B$.   See, for example,
Section 1.1.10, item 7, page 11 of [\cite{GouldenJackson}], where
the identity is attributed to Jacobi.

\BeginProofOf{Theorem \ref{theorem3}}  We start with the expected
determinant,
 \begin{align*}
&~~~~ E(\det(G_{n})) \\
 &= \left[ \frac{x^{n}}{n!} \right] \exp{\left\lbrace\frac{t_{1}x}{1}-\frac{t_{2}x^{2}}{2}+\frac{t_{3}x^{3}}{3}-
\dots\right\rbrace} .\\
&= \left[ \frac{x^{n}}{n!} \right] \exp{ \left\lbrace \trace\left(
\frac{M^{1}x}{1}-\frac{M^{2}x^{2}}{2}+\frac{M^{3}x^{3}}{3}- \dots
\right) \right\rbrace} \\
&= \left[ \frac{x^{n}}{n!} \right] \det{ \left( \exp{\left\lbrace
 \frac{M^{1}x}{1}-\frac{M^{2}x^{2}}{2}+\frac{M^{3}x^{3}}{3}-
\dots \right\rbrace} \right)}  \\
&= \left[ \frac{x^{n}}{n!} \right] \det{ \left( \exp \left\lbrace
\log{(I+xM)} \right\rbrace \right)} \\
&= \left[ \frac{x^{n}}{n!} \right] \det{(I+xM)}  \\ 
&= \left[ \frac{x^{n}}{n!} \right]  (-x)^{t} \cdot \det{(\lambda
I-M) \vert_{\lambda=-\frac{1}{x}}} \\
&= \left[ \frac{x^{n}}{n!} \right]  (-x)^{t} \left( \lambda^{t}
- c_{1}\lambda^{t-1} + \cdots \right) \vert_{\lambda=-\frac{1}{x}}
\\
&= \left[ \frac{x^{n}}{n!} \right]  (-x)^{t} \left( \left(-\frac{1}{x}\right)^{t}
- c_{1}\left(-\frac{1}{x}\right)^{t-1} + \cdots \right) \\
 &= \left[ \frac{x^{n}}{n!} \right] \left(1 + c_{1}x + c_{2}x^{2}
+ c_{3}x^{3} + \cdots \right)	 \\
 &= n!c_{n} \; \; \nonumber .
\end{align*}
The first equality comes from (\ref{eqgfdet}), the third from
(\ref{eqexptrace}), and the rest are straightforward
manipulations.  The proof for the expected permanent is similar:
\begin{align*}
 &~~~~E(\perm(G_{n})) \\
 &= \left[ \frac{x^{n}}{n!} \right] \exp{\left\lbrace\frac{t_{1}x}{1}+\frac{t_{2}x^{2}}{2}+\frac{t_{3}x^{3}}{3}+
\dots\right\rbrace} \\
&= \left[ \frac{x^{n}}{n!} \right] \exp{ \left\lbrace \trace\left(
\frac{M^{1}x}{1}+\frac{M^{2}x^{2}}{2}+\frac{M^{3}x^{3}}{3}+ \dots
\right) \right\rbrace} \\
&= \left[ \frac{x^{n}}{n!} \right] \det{ \left( \exp{\left\lbrace
 \frac{M^{1}x}{1}+\frac{M^{2}x^{2}}{2}+\frac{M^{3}x^{3}}{3}+
\dots \right\rbrace} \right)}  \\
&= \left[ \frac{x^{n}}{n!} \right] \det{ \left( \exp \left\lbrace
\log{(I-xM)^{-1}} \right\rbrace \right)} \\
&= \left[ \frac{x^{n}}{n!} \right] \det{(I-xM)^{-1}}  \\
&= \left[ \frac{x^{n}}{n!} \right] \frac{1}{\det{(I-xM)}} \\
&= \left[ \frac{x^{n}}{n!} \right]  \dfrac{1}{(x)^{t} \cdot \det{(\lambda
I-M)} \vert_{\lambda=\frac{1}{x}}} \\
&= \left[ \frac{x^{n}}{n!} \right]  \frac{1}{(x)^{t} \left( \left(\frac{1}{x}\right)^{t}
- c_{1}\left(\frac{1}{x}\right)^{t-1} + \cdots \right)} \\
&= \left[ \frac{x^{n}}{n!} \right]  \frac{1}{
             1 - c_{1}x + c_2x^2 - c_3x^3 + \cdots
                                             }
\end{align*}
This time the first equality follows from (\ref{eqgfperm}), the third
again is
from (\ref{eqexptrace}), and, as was the case with the determinant,
the rest are straightforward manipulations.  
\EndProof

\noindent{\bf Remark.}  The function
$$
\zeta_G(u) =  \frac{1}{\det(I-Tu)},
$$
where $T$ is Hashimoto's edge adjacency operator,
is called the {\it Ihara zeta-function of the graph}
$G$, see \cite{Wiki}, \cite{Terras}.

\vskip 10pt

\BeginProofOf{Theorem \ref{theorem4}}
For simplicity, let us assume the probability distribution on 
vectors $w$ is discrete;
say, $v_{1}$, $v_2$, $\dots$ with probabilities $p(v_{1})$,
$p(v_{2})$, $\dots$.
In order to obtain a term $\lambda^{n-t}$ in the expansion of
$\det(\lambda I - A^{T}A)$, we choose the $n-t$ $\lambda$'s from the
main diagonal, and then expand the remaining principal submatrix of
size $t$.  Since the remaining submatrix is that of $-A^TA$, we obtain
the sign-adjustment $(-1)^t$ and find
\begin{equation*}
 	b_{i} = \sum_{\substack{i \times i \\ \textnormal{principal}
\\ \textnormal{submatrices } \alpha}}\det(\alpha) \; \; .
\end{equation*}
Since expectation is a linear operator, the expected
value of the $i$'th coefficient of the characteristic polynomial
of $A^{T}A$ is 
\begin{align*}
 E(b_{i})&= E\left(\sum_{\substack{i \times i \\ \textnormal{principle}
\\ \textnormal{submatrices } \alpha}} \det(\alpha)\right) \\
&= \binom{n}{i} \sum_{(v_{1},\dots,v_{i})} \det{\left( \begin{array}{cccc}
(v_{1}\cdot v_{1}) & (v_{1}\cdot v_{2}) &  \cdots & (v_{1}\cdot
v_{i}) \\ (v_{2}\cdot v_{1}) & (v_{2}\cdot v_{2}) &  \cdots &
(v_{2}\cdot v_{i}) \\ \vdots & \vdots  & \vdots & \vdots \\ (v_{i}\cdot
v_{1}) & (v_{i}\cdot v_{2}) &  \cdots & (v_{i}\cdot v_{i})  \end{array}\right)}
p(v_{1}) \dots p(v_{i})  \\
&= \binom{n}{i} E(\det(G_{i})) \; \; .
\end{align*}
This proves the first assertion of the theorem, and the second assertion
regarding the permanental polynomial is demonstrated in a similar manner.
\EndProof

\section{Experimental results}
We wrote a Matlab program to compare the expected characteristic
and permanental polynomials given by Theorem~\ref{theorem4}
to those of randomly sampled matrices of various
sizes.  We computed all permanents using a programmatic link with
Maple via the Maple Toolbox for Matlab, and all characteristic
polynomials using the Matlab command \textit{poly}.

There are, of course, infinitely many different distributions which
might underly the vectors $w$.   We chose to use one based
on {\it sample counts}.
This is an easily understood distribution, and of interest for
possible applications.   The idea is to assume a set $X$ of size
$t$, $\{x_1,\dots,x_t\}$, with probabilities $p_1,\dots, p_t$,
where $\sum_{i=1}^tp_i=1$.  We take a sample, with replacement,
from $X$ of size $\ell$, and record $w_i$ as the number of times
that the element $x_i$ is chosen.  Then, each of our $t$-tall
vectors $w$ is integral, satisfies $\sum_{i=1}^tw_i=\ell$,
and the distribution on these vectors is the familiar multinomial
distribution:
\begin{equation*}
 \textnormal{Prob}\left\lbrace \left[ \begin{array}{c} w_{1}
\\  w_{2} \\ \vdots \\ w_{t} \end{array} \right] =\left[ \begin{array}{c}
b_{1} \\  b_{2} \\ \vdots \\ b_{t} \end{array} \right] \right\rbrace
= {\binom{l}{b_{1} \dots b_{t}}} p_{1}^{b_{1}} \cdots p_{t}^{b_{t}}
\; \; .
\end{equation*}
The corresponding matrix $M$ of second moments is found to be
\begin{equation*}
 M_{ij} = E(w_{i}w_{j}) =
\begin{cases}
		\ell(\ell-1)p_{i}^{2} + \ell p_{i} & \textnormal{if } \; i=j\\
                \ell(\ell-1)p_{i}p_{j} & \textnormal{if } \; i\neq j \; .
  	         \end{cases}
\end{equation*}
We have derived the matrix $M$ for several other scenarios which seem
natural for applications, but do not report any of these results in the
present paper, with one exception.  Namely, suppose that the
random vectors $w$ are generated as counts, much as above, except
the sample size $\ell$ is also a random quantity.  That is,
the $w$ come about by a compound process.  If we assume the sample
size $\ell$ to be given by a distribution $\prob(\ell)$, then
$$
M_{ij} = \begin{cases}
			\sum_{\ell}\left( \ell(\ell-1)p_{i}^{2}+\ell p_{i}\right) \prob(\ell),
& i=j\\
				\sum_{\ell}\left(\ell(\ell-1)p_{i}p_{j}\right) \prob(\ell)
, &  i \neq j  \\ 
		\end{cases}\; \;.
$$
In our experiments, we used the above multinomial distribution
to generate random vectors with $\ell=10$, $t=4$
and $p= \left[ 3/8, 1/4, 1/4, 1/8 \right]$.   

\begin{figure}[htbp]
 \centering
 \includegraphics[height=5.5in, width=6.5in]{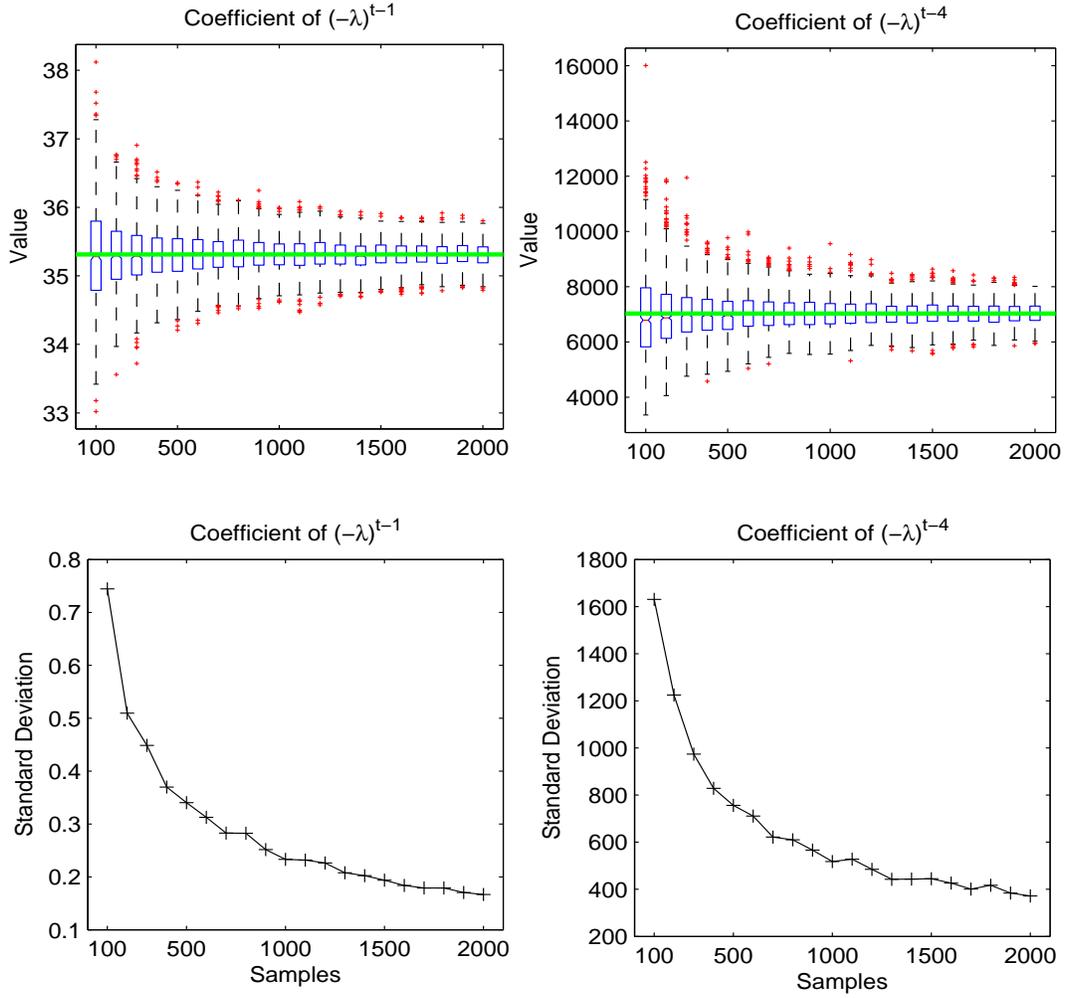}
\caption{Experimental results of the coefficients of \textbf{characteristic} polynomials
of $1000$ matrices for particular sample sizes (columns) with
the model $l=10$, $t=4$ and $p= \left[ 3/8, 1/4, 1/4, 1/8 \right]$.
 To normalize the results, each coefficient of $(-\lambda)^{t-i}$
was divided by  $\binom{n}{i}$ before plotting.  Green line represents
$E(\det(G_{i}))$ as computed by the recursion (\ref{eqdetrecursion}).}
\label{fig:detcoeffs}
 \end{figure}

We know, theoretically, that $a_{i}=0$ for $i > t$.  However,
we computed these for confirmation.  Thus, we computed $t_{i}$
for $i$ up to $7$.
\begin{center}
\begin{tabular}{l|c}
$i$ & $t_{i}$ \\ \hline 
1 & 565/16 \\ 
2 & 210825/256 \\ 
3 & 93917125/4096 \\ 
4 & 42581180625/65536 \\ 
5 & 19338382478125/1048576 \\ 
6 & 8784040432265625/16777216 \\ 
7 & 3990026079685703125/268435456
\end{tabular}
\end{center}
Then, from recursion (\ref{eqdetrecursion}), we compute $a_i$.
We also generated $1000$ matrices $G_{4000}=A^{T}A$ at random
(by sampling $4000$ random vectors from the counting distribution
described above to form $A$) and computed the mean and standard deviation of
the appropriate coefficient
in  $\det(\lambda I - G_{4000})$, divided by a binomial coefficient
as given in equation (\ref{eqcharcoefexpect}).  The exact values
and the sample values were then compared.  Here are the results:

\begin{center}
\begin{tabular}{l|ccc}
$i$ & $a_{i}$, recursion & $a_{i}$, sample & std dev \\ \hline
1 & 35.31 & 35.32 & 0.12 \\ 
2 &  423.44 & 423.59 & 5.91 \\ 
3 &  2648.44 & 2649.99 & 67.54 \\ 
4 &  7031.25 & 7035.47 & 263.93 \\ 
5 &  0.0 & 0.0 & 0.0 \\ 
6 &  0.0 & 0.0 & 0.0
\end{tabular} 
\end{center} 
Furthermore, boxplots of the sample coefficient distributions
as the sample size increases are depicted in Figure \ref{fig:detcoeffs}.
These plots show that, for this distribution, the standard deviation
decreases according to a power law.

\vskip 10pt

\begin{figure}[htbp]
 \centering
 \includegraphics[height=5.5in, width=6.5in]{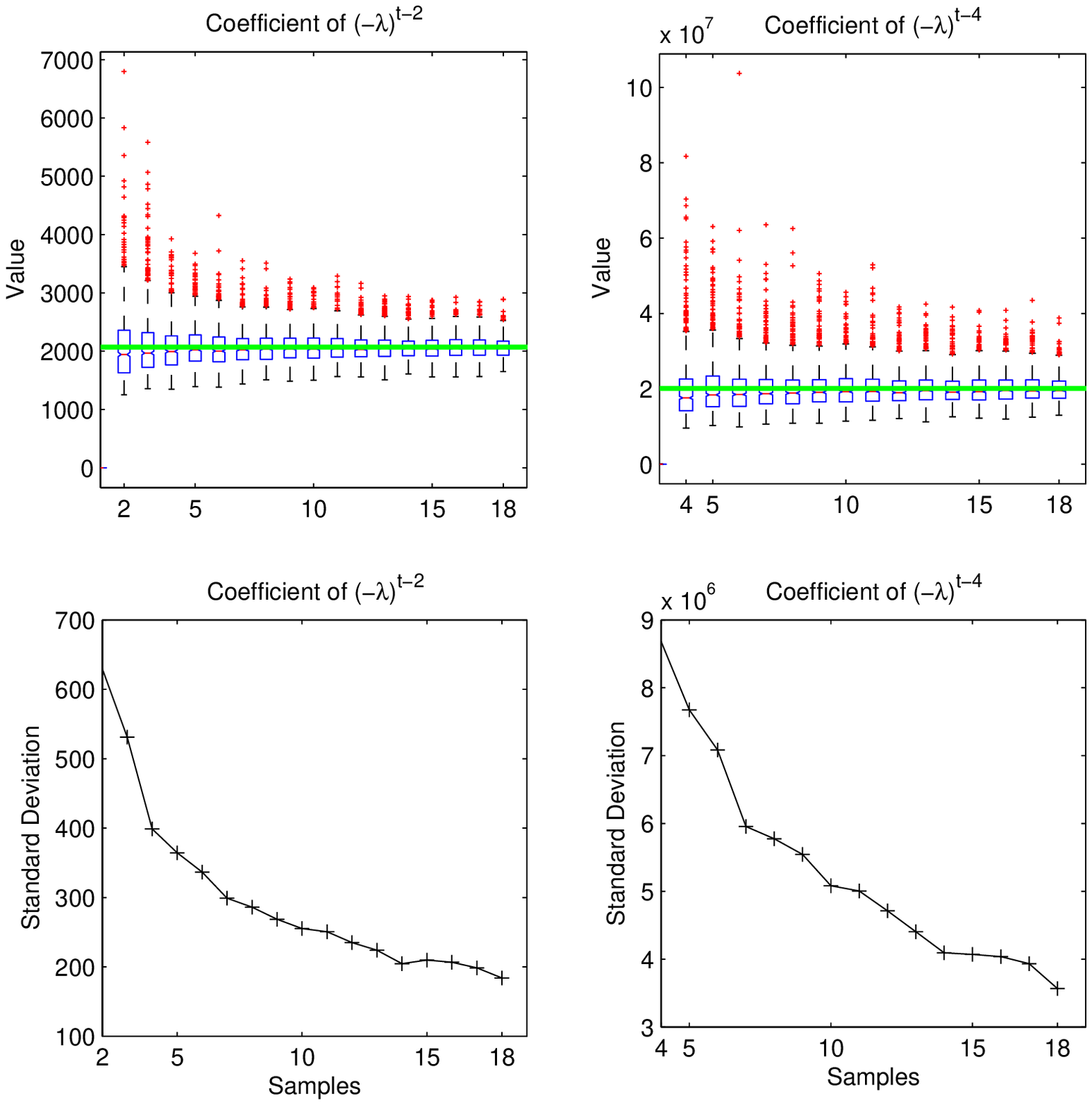}
\caption{Experimental results of the \textbf{permanental} coefficients
from $1000$ matrices for particular sample sizes (columns) with
the model $l=10$, $t=4$ and $p= \left[ 3/8, 1/4, 1/4, 1/8 \right]$.
 To normalize the results, each coefficient of $(-\lambda)^{t-i}$
was divided by  $\binom{n}{i}$ before plotting.  Green line represents
$E(\perm(G_{i}))$ as computed by the recursion (\ref{eqpermrecursion}).}
\label{fig:permcoeffs}
\end{figure}

\noindent
Using the same $t_{n}$ as above, we used the permanental coefficient
recursion (\ref{eqpermrecursion}) to compute the exact values of
$p_{n}=E(\perm(G_{n}))$.
We also computed the coefficients of permanental polynomials of $1000$ random
gram matrices $G_{18}=A^{T}A$ (created by sampling $18$ random
vectors from the distribution) and subsequently compared the
sampled results with those provided by the recursion.  
We computed the mean and standard deviation after division by
the binomial coefficient as given in (\ref{eqpermcoefexpect}).
The exact and sampled values were then compared.  Due to the
intractable computational complexity of computing the exact permanent,
we were computationally limited to computing only matrices with
$18$ randomly sampled columns.  Here are the results (see
Figure \ref{fig:permcoeffs} for boxplots):
\begin{center}
\begin{tabular}{l|ccc}
$i$ & $p_{i}$, recursion & $p_{i}$, sample & std dev \\ \hline
1 & 35.31 & 35.27 & 1.921 \\ 
2 &  2070.51 & 2071.54 & 206.63 \\ 
3 &  177134.95 & 177679.22 & 26299.36 \\ 
4 &  20126988.14 & 20245985.44 & 4037955.48 \\ 
5 &  2857210195.90 & 2882490271.36 & 729446452.08 \\ 
6 &  486697830067.95 & 492457249647.11 & 151888811726.52
\end{tabular}
\end{center}

\section{Connection to the cycle index polynomial}

For a permutation $\sigma\in S_n$ belonging to the symmetric group of
order $n$, let $N_i(\sigma)$ denote the number of cycles in $\sigma$
of size $i$.
Let $X_1,X_2,\dots$ be a countably infinite sequence of variables, and
define the polynomial $P_n(X_1,\dots,X_n)$ by
$$
P_n(X_1,\dots,X_n) = \sum_{\sigma\in S_n} \prod_{i=1}^n X_i^{N_i(\sigma)}.
$$
Then the quotient $P_n(X_1,\dots,X_n)/n!$ is called the {\it cycle index
polynomial of the symmetric group}.   The generating function identity
$$
\sum_{n=0}^{\infty} P_n(X_1,\dots,X_n) \frac{u^n}{n!} ~~=~~
\exp\left(
               \sum_{i=1}^{\infty} \frac{X_iu^i}{i}
   \right)
$$
was observed in \cite{BC}.  The latter paper was devoted to proving
that assigning nonnegative real values to the variables $X_i$ subject to
certain inequalities would result in the real values $P_n(X_1,\dots,X_n)$
satisfying similar inequalities.  Coincidentally, the pairs of
quantities $(-1)^{n-1}t_n,a_n$ and $t_n, p_n$
studied in this paper satisfy identical generating function identities.
In particular, the sequence of expected permanents $p_n=E(\perm(A^TA))$
are hereby identified as evaluations of the cycle index polynomials at
certain weights $t_i$.

\section{Summary and conclusion}\label{lblconclusion}

We have introduced the notion of a random Gram matrix, and
provided theory enabling the efficient computation of the expected
determinant and expected permanent of it.
The random Gram matrix consists of dot products of vectors taken
from various distributions.  We further
proved generating function identities and recursions relating
these expectations to the traces of powers of a second moment matrix.
The expected coefficients of the characteristic and
permanental polynomials have also been studied, with some numerical
experiments checking on the theory.
Some of the formulas found are the same as
those studied in earlier work in an entirely
different context \cite{BC}.

We have observed
empirically that as the number of columns
in the sample matrix $A$ increases, the standard deviation of the
normalized expected coefficients of the determinantal
and permanental polynomials decreases according to a power law. 
Although the empirical data presented in this paper was limited
to the multinomial counting model, the theoretical relationships
between the different quantities
remain no matter which representation is used.  
In future work, the theoretical rate of convergence should
be formulated according to the representation and probability
model used to generate the matrix $A$ (e.g. trivially, when A=0, 
the truth converges immediately to the expected value).  

Can the probabilistic results presented in this work be of any
help in managing the complexity of computing
the permanent?  Already, \cite{Barvinok}, there is a polynomial
time algorithm for computing the permanent of an
$n\times n$ matrix of rank $t$, $t$ being fixed.  One way for
our probabilistic methods to impact complexity considerations
would be via finding a distribution on $t$-vectors ($t$ small) such that
a given $n\times n$ permanent per$(H)$ is equal to or
well approximated by the expected value of per$(A^T A)$.
We have no ideas in this direction.

It is hoped
that the theoretical observations we have made
will prove useful in processing and comparing large amounts of
numerical data, such as those algorithms that use permanental
polynomials of large chemical graphs \cite{Cash,Cash2}.  Moreover,
the combinatorial relationships between traces of matrix powers,
characteristic coefficients, expected permanents, and expected
determinants will help us better understand how to use these
quantities, create bounds for them, and illuminate what has made
them so especially useful in applied numerical science.  


\end{document}